\documentclass{article}%
\usepackage{typearea}
\usepackage{enumerate}
\usepackage{amsmath}%
\usepackage{amsfonts}%
\usepackage{amssymb}%
\usepackage[center]{subfigure}
\usepackage{todonotes}
\usepackage[font=small,format=plain,labelfont=bf,up,textfont=it,up]{caption}
\usepackage{enumitem}
\usepackage{graphicx}
\usepackage[pdftex]{geometry}
\usepackage{mathtools}
\usepackage{color}%

\providecommand{\U}[1]{\protect\rule{.1in}{.1in}}
\newtheorem{theorem}{Theorem}

\newtheorem{corollary}[theorem]{Corollary}

\newtheorem{lemma}[theorem]{Lemma}

\newtheorem{remark}[theorem]{Remark}

\newenvironment{proof}[1][Proof]{\noindent\textbf{#1.} }{\ \rule{0.5em}{0.5em}}
\newcommand{\BIGOP}[1]{\mathop{\mathchoice{\raise-0.22em\hbox{\huge
$#1$}} {\raise-0.05em\hbox{\Large $#1$}}{\hbox{\large $#1$}}{#1}}}

\DeclareMathOperator*{\argmin}{arg\,min}

\begin{document}

\title{Numerical approximation of Poisson problems in long domains}
\author{M. Chipot\thanks{Institute for Mathematics, University of Z\"{u}rich,
Winterthurerstr. 190, CH-8057 Z\"{u}rich \texttt{m.m.chipot@math.uzh.ch}}
\and W. Hackbusch\thanks{Max Planck Institute for Mathematics in the Sciences,
Inselstr. 22, D-04103 Leipzig, \texttt{wh@mis.mpg.de}}
\and S. Sauter\thanks{Institute for Mathematics, University of Z\"{u}rich,
Winterthurerstr. 190, CH-8057 Z\"{u}rich \texttt{stas@math.uzh.ch}}
\and A. Veit\thanks{Institute for Mathematics, University of Z\"{u}rich,
Winterthurerstr. 190, CH-8057 Z\"{u}rich \texttt{alexander.veit@math.uzh.ch}}}
\date{}
\maketitle

\begin{abstract}
In this paper, we consider the Poisson equation on a \textquotedblleft
long\textquotedblright\ domain which is the Cartesian product of a
one-dimensional long interval with a $\left(  d-1\right)  $-dimensional
domain. The right-hand side is assumed to have a rank-1 tensor structure. We
will present and compare methods to construct approximations of the solution which have
tensor structure and the computational effort is governed by only solving
elliptic problems on lower-dimensional domains. A zero-th order tensor
approximation is derived by using tools from asymptotic analysis (method 1).
The resulting approximation is an elementary tensor and, hence has a fixed
error\ which turns out to be very close to the best possible approximation of
zero-th order. This approximation can be used as a starting guess for the
derivation of higher-order tensor approximations by an
alternating-least-squares (ALS) type method (method 2). Numerical experiments
show that the ALS is converging towards the exact solution.

Method 3 is based on the derivation of a tensor approximation via exponential
sums applied to discretised differential operators and their inverses. It can
be proved that this method converges exponentially with respect to the tensor
rank. We present numerical experiments which compare the performance and
sensitivity of these three methods.



\end{abstract}

\section{Introduction}

In this paper, we consider elliptic partial differential equations on domains
which are the Cartesian product of a \textquotedblleft long\textquotedblright%
\ interval $I_{\ell}=\left(  -\ell,\ell\right)  $ with a $\left(  d-1\right)
$-dimensional domain $\omega$, the \textit{cross section} - a typical
application is the modelling of a flow in long cylinders. As a model problem
we consider the Poisson equation with homogeneous Dirichlet boundary
conditions and a right-hand side which is an \textit{elementary tensor}; i.e.,
the product of a univariate function (on the long interval) and a $\left(
d-1\right)  $-variate function on the cross section. Such problems have been
studied by using asymptotic analysis, see., e.g., \cite{chipot2016_asym}. Our
first approximation (method 1) is based on this technique and approximates the
solution by an elementary tensor where the function on the cross section is
the solution of a Poisson-problem on the cross section and the corresponding
univariate function is determined afterwards as the best approximation in the
Sobolev space $H_{0}^{1}$ on the long interval. In Lemma
\ref{convergenceResultM1} below, it is shown that this approximation converges
exponentially with respect to the length of the cylinder for any subdomain
$I_{\ell_{0}}\times\omega$ for fixed $\ell_{0}<\ell$. However, for fixed
$\ell$ this is a one-term approximation with a fixed error.

Method 2 uses the result of method 1 as the initial guess for an iterative
procedure which is an \textit{alternating least squares} (ALS) type method.
Recursively, one assumes that a rank-$k$ tensor approximation of the solution
has already been derived and then starts an iteration to compute the $k+1$
term: a) one chooses an univariate function on $I_{\ell}$ as an initial guess
for this iteration and determines the function on the cross section as the
best approximation in $H_{0}^{1}$ of the cross-section. In step b) the
iteration is flipped and one fixes the new function on the cross section and
determines the corresponding best approximation in $H_{0}^{1}$ of the
interval. Steps a) and b) are iterated until a stopping criterion is reached
and this gives the $k+1$ term in the tensor approximation. In the literature this approach is also known as \emph{Proper Generalized Decomposition} (PGD) (\cite{chinesta2013proper}).  We have performed
numerical experiments which are reported in Section \ref{NumericalExperiments}
which show that this method leads to a convergent approximation also for fixed
$\ell$ as the tensor rank of the approximation increases. However, it turns
out that this method is quite sensitive and requires that the inner iteration
a), b) leads to an accurate approximation of the $\left(  k+1\right)  $ term
in order to ensure that the outer iteration is converging. Furthermore, the
numerical experiments that we have performed indicate that the convergence
speed can slow down as the number of outer iterations increases. Thus, this
method is best suited when a medium approximation accuracy of the Poisson
problem is required. 


Method 3 is based on a different approach which employs numerical tensor
calculus (see \cite{hackbusch2012tensor}). First one defines an exponential
sum approximation of the function $1/x$. Since the differential operator
$-\Delta$ is of tensor form, the exponential sum, applied to the inverse of a
discretisation of the Laplacian by a matrix which must preserve the tensor
format, directly leads to a tensor approximation of the solution $u$. We
emphasize that the explicit computation of the inverse of the discretisation
matrix can be avoided by using the \textit{hierarchical format} for their
representation (see \cite{hackbusch2015hierarchical}). An advantage of this
method is that a full theory is available which applies to our application and
allows us to choose the tensor rank via an a priori error estimate. It also
can be shown that the tensor approximation converges exponentially with
respect to the tensor rank (see \cite{hackbusch2012tensor}).

The goal of this paper is to compare three different approaches for the numerical approximation of Poisson problems in domains of the form $I_\ell\times\omega$ and to assess their performance with respect to the length $\ell$ via numerical experiments. These methods exhibit different computational complexities and our results can be used to determine a suitable method given a desired accuracy range. For an in-depth theoretical analysis of the presented methods we refer to the existing literature.

The paper is structured as follows. In Section \ref{ProblemFormulation} we
formulate the problem on the long product domain and introduce the assumptions
on the \textit{tensor format} of the right-hand side. The three different
methods for constructing a tensor approximation of the solution are presented
in Section \ref{NumericlaApproximation}. The results of numerical experiments
are presented in Section \ref{NumericalExperiments} where the convergence and
sensitivity of the different methods is investigated and compared. For the
experiments we consider first the case that the cross section is the
one-dimensional unit interval and then the more complicated case that the
cross section is an L-shaped polygonal domain. Finally, in the concluding
section we summarize the results and give an outlook.

\section{Setting\label{ProblemFormulation}}

Let $\omega$ be an open, bounded and connected Lipschitz domain in
$\mathbb{R}^{n-1},n\geq1$. In the following we consider Poisson problems on
domains of the form
\[
\Omega_{\ell}:=I_{\ell}\times\omega\quad\text{with\quad}I_{\ell}:=\left(
-\ell,\ell\right)  ,
\]
where $\ell$ is large. We are interested in Dirichlet
boundary value problems of the form
\begin{align}
-\Delta u_{\ell} &  =F\quad\text{ in }\Omega_{\ell},\nonumber\\
u_{\ell} &  =0\quad\text{ on }\partial\Omega_{\ell}\label{originalDBVP}%
\end{align}
with weak formulation%
\[%
\begin{cases}
\text{find }u_{\ell}\in H_{0}^{1}(\Omega_{\ell})\text{ s.t. }\\
(\nabla u_{\ell},\nabla v)_{L_{2}(\Omega_{\ell})}=(F,v)_{L_{2}(\Omega_{\ell}%
)}\quad\forall v\in H_{0}^{1}(\Omega_{\ell}).
\end{cases}
\]
Specifically we are interested in right-hand sides $f$ which have a tensor
structure of the form $F=1\otimes f$ or more generally $F=\sum_{k=0}^{n}%
g_{k}\otimes f_{k}$, where $g_{k}$ is a univariate function and $f$, $f_{k}$
are functions which depend only on the $\left(  d-1\right)  $-dimensional
variable $x^{\prime}\in\omega$. Here, we use the standard tensor notation
$\left(  g\otimes f\right)  \left(  x\right)  =g\left(  x_{1}\right)  f\left(
x^{\prime}\right)  $ with $x^{\prime}=\left(  x_{k}\right)  _{k=2}^{d}$. In
this paper, we will present and compare methods to approximate $u_{\ell}$ in
tensor form.\medskip

\noindent We consider a right-hand side of the form%
\begin{equation}
F=1\otimes f\qquad\text{for some }f\in L^{2}\left(  \omega\right)
\label{defFtensor1}%
\end{equation}
and derive a first approximation of $u_{\ell}$ as the solution of the $\left(
n-1\right)  $-dimensional problem on $\omega$:
\begin{align}
-\Delta^{\prime}u_{\infty}(x^{\prime})  &  =f(x^{\prime})\quad\text{ in
}\omega,\nonumber\\
u_{\infty}  &  =0\quad\text{ on }\partial\omega\label{reducedDBVP}%
\end{align}
with weak form%
\[%
\begin{cases}
\text{find }u_{\infty}\in H_{0}^{1}(\omega)\text{ s.t. }\\
(\nabla^{\prime}u_{\infty},\nabla^{\prime}v)_{L_{2}(\omega)}=(f,v)_{L_{2}%
(\omega)}\quad\forall v\in H_{0}^{1}(\omega).
\end{cases}
\]

\section{Numerical Approximation\label{NumericlaApproximation}}

In this section we derive three different methods to approximate problem
(\ref{originalDBVP}). In all three methods we exploit the special structure of
the domain $\Omega_{\ell}$ and the right-hand side $F$. Our goal is to reduce
the original $n$-dimensional problem on $\Omega_{\ell}$ to one or more
$\left(  n-1\right)  $-dimensional problems on $\omega$. Compared to standard
methods like finite elements methods or finite difference methods, which solve
the equations on $\Omega_{\ell}$, this strategy can significantly reduce the
computational cost since $\ell$ is considered large and the discretisation in
the $x_{1}$ direction can be avoided.

\subsection{Method 1: A one-term approximation based on an asymptotic analysis
of problem (\ref{originalDBVP})}

Although the right-hand side $F$ in (\ref{originalDBVP}) is independent of
$x_{1}$, it is easy to see that this is not the case for the solution
$u_{\ell}$, i.e., due to the homogeneous Dirichlet boundary conditions it is
clear that $u_{\ell}$ depends on $x_{1}$. However, if $\ell$ is large one can
expect that $u_{\ell}$ is approximately constant with respect to $x_{1}$ in a
subdomain $\Omega_{\ell_{0}}$, where $0<\ell_{0}\ll\ell$ and thus converges
locally to a function independent of $x_{1}$ for $\ell\rightarrow\infty$. The
asymptotic behaviour of the solution $u_{\ell}$ when $\ell\rightarrow\infty$
has been investigated in \cite{chipot2009elliptic}. It can be shown that
\[
u_{\ell}\longrightarrow1\otimes u_{\infty}\qquad\text{in }\Omega_{\ell_{0}},
\]
where $u_{\infty}$ is the solution of (\ref{reducedDBVP}), with an exponential
rate of convergence. More precisely, the following theorem holds:

\begin{theorem}
\label{u_inf_estimate}There exist constants $c,\alpha>0$ independent of $\ell$
s.t.
\[
\int_{\Omega_{\ell/2}}|\nabla(u_{\ell}-1\otimes u_{\infty})|^{2}dx\leq
c\operatorname{e}^{-\alpha\ell}\Vert f\Vert_{2,\omega}^{2},
\]
where $\Vert\cdot\Vert_{2,\omega}$ refers to the $L^{2}(\omega)$-norm.
\end{theorem}

For a proof we refer to \cite[Theorem 6.6]{chipot2009elliptic}.

Theorem \ref{u_inf_estimate} shows that $1\otimes u_{\infty}$ is a good
approximation of $u_{\ell}$ in $\Omega_{\ell/2}$ when $\ell$ is large. This
motivates to seek approximations of $u_{\ell}$ in $\Omega_{\ell}$ which are of
the form
\[
u_{\ell}\approx u_{\ell}^{M_{1}}:=\psi_{\ell}\otimes u_{\infty},
\]
where $\psi_{\ell}\in H_{0}^{1}(-\ell,\ell)$. Here, we choose $\psi_{\ell}$ to
be the solution of the following best approximation problem: Given $u_{\ell
}\in H_{0}^{1}(\Omega_{\ell})$ and $u_{\infty}\in H_{0}^{1}(\omega)$, find
$\psi\in H_{0}^{1}(-\ell,\ell)$ s.t.
\begin{equation}
\Vert\nabla\left(  u_{\ell}-\psi\otimes u_{\infty}\right)  \Vert_{2}=
\inf_{\theta\in H_{0}^{1}(-\ell,\ell)}\Vert\nabla\left(  u_{\ell}%
-\theta\otimes u_{\infty}\right)  \Vert_{2}. \label{bestApproximation1}%
\end{equation}
In order to solve problem (\ref{bestApproximation1}) we define the functional
\[
J(u_{\ell},u_{\infty})(\theta):=\left\Vert \nabla\left(  u_{\ell}%
-\theta\otimes u_{\infty}\right)  \right\Vert _{2}^{2}%
\]
and consider the variational problem of minimizing it with respect to
$\theta\in H_{0}^{1}\left(  -\ell,\ell\right)  $.

A simple computation shows that this is equivalent to finding $\tilde{\theta
}\in H_{0}^{1}\left(  I_{\ell}\right)  $ such that
\[
\left.
\begin{array}
[c]{cl}
& \left(  \nabla\left(  \theta\otimes u_{\infty}\right)  ,\nabla\left(
\tilde{\theta}\otimes u_{\infty}\right)  \right)  _{2}=\left(  \nabla u_{\ell
},\nabla\left(  \tilde{\theta}\otimes u_{\infty}\right)  \right)  _{2}\\
\iff & \left(  \left(
\begin{array}
[c]{c}%
\theta^{\prime}\otimes u_{\infty}\\
\theta\otimes\nabla^{\prime}u_{\infty}%
\end{array}
\right)  ,\left(
\begin{array}
[c]{c}%
\tilde{\theta}^{\prime}\otimes u_{\infty}\\
\tilde{\theta}\otimes\nabla^{\prime}u_{\infty}%
\end{array}
\right)  \right)  _{2}=\left(  -\Delta u_{\ell},\tilde{\theta}\otimes
u_{\infty}\right)  _{2}\\
\iff & \alpha_{\infty}^{2}\left(  \theta^{\prime},\tilde{\theta}^{\prime
}\right)  _{2,I_{\ell}}+\beta_{\infty}^{2}\left(  \theta,\tilde{\theta
}\right)  _{2,I_{\ell}}=\left(  1\otimes f,\tilde{\theta}\otimes u_{\infty
}\right)  _{2}\\
\iff & \alpha_{\infty}^{2}\left(  \theta^{\prime},\tilde{\theta}^{\prime
}\right)  _{2,I_{\ell}}+\beta_{\infty}^{2}\left(  \theta,\tilde{\theta
}\right)  _{2,I_{\ell}}=\left(  1\otimes\left(  -\Delta^{\prime}u_{\infty
}\right)  ,\tilde{\theta}\otimes u_{\infty}\right)  _{2}\\
\iff & \alpha_{\infty}^{2}\left(  \theta^{\prime},\tilde{\theta}^{\prime
}\right)  _{2,I_{\ell}}+\beta_{\infty}^{2}\left(  \theta,\tilde{\theta
}\right)  _{2,I_{\ell}}=\beta_{\infty}^2\int_{I_{\ell}}\tilde{\theta}.
\end{array}
\right\}  \;\forall\tilde{\theta}\in H_{0}^{1}\left(  I_{\ell}\right)
\]
with%
\[
\alpha_{\infty}:=\left\Vert u_{\infty}\right\Vert _{2,\omega},\quad
\beta_{\infty}:=\left\Vert \nabla^{\prime}u_{\infty}\right\Vert _{2,\omega}.
\]
The strong form of the resulting equation is%
\begin{align*}
-\alpha_{\infty}^{2}\theta^{\prime\prime}+\beta_{\infty}^{2}\theta &
=\beta_{\infty}^{2}\quad\text{ on }(-\ell,\ell),\\
\theta(-\ell)=\theta(\ell)  &  =0.
\end{align*}
The solution of this one-dimensional boundary value problem is given by
\[
\theta\left(  x_{1}\right)  :=1-\frac{\cosh\left(  \frac{\beta_{\infty}%
}{\alpha_{\infty}}x_{1}\right)  }{\cosh\left(  \frac{\beta_{\infty}}%
{\alpha_{\infty}}\ell\right)  }.
\]
This shows that an approximation of our original problem (\ref{originalDBVP})
is given by
\begin{equation}
u_{\ell}^{M_{1}}:=\psi_{\ell}\left(  \lambda_{\infty},\cdot\right)  \otimes
u_{\infty},\quad\text{\ with\ }\psi_{\ell}(a,x_{1}):=1-\frac{\cosh\left(
ax_{1}\right)  }{\cosh\left(  a\ell\right)  } \label{approximationM1}%
\end{equation}
and%
\begin{equation}
\lambda_{\infty}=\frac{\beta_{\infty}}{\alpha_{\infty}}=\frac{\sqrt{\left(
f,u_{\infty}\right)  _{2,\omega}}}{\alpha_{\infty}}. \label{lambdainf}%
\end{equation}
Note that $\psi_{\ell}\left(  a,\cdot\right)  $ satisfies%
\begin{equation}
-\psi_{\ell}^{\prime\prime}\left(  a,\cdot\right)  +a^{2}\psi_{\ell}\left(
a,\cdot\right)  =a^{2}\quad\text{and\quad}\psi_{\ell}\left(  \pm\ell\right)
=0. \label{dglpsi}%
\end{equation}
In Section \ref{NumericalExperiments} we report on various numerical
experiments that show the approximation properties of this rather simple
one-term approximation.

\begin{figure}[tbh]
\centering
\includegraphics[width=0.95\textwidth]{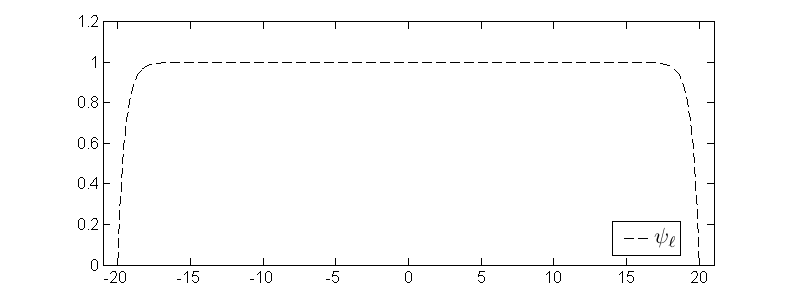} \caption{Plot of $\psi
_{\ell}\left(  \lambda_{\infty},\cdot\right)  $ for $\ell=20$ and
$\lambda_{\infty}=2$}%
\label{plotGamma}%
\end{figure}

Figure \ref{plotGamma} shows a plot of $\psi_{\ell}\left(  \lambda_{\infty
},\cdot\right)  $ for $\ell=20$ and $\lambda_{\infty}=2$. Since $\psi_{\ell}$
approaches $1$ with an exponential rate as $x_{1}$ moves away from $\pm\ell$
towards the origin, an analogous result to Theorem \ref{u_inf_estimate} can be
shown for $u_{\ell}^{M_{1}}$.

\begin{lemma}
\label{convergenceResultM1} There exist constants $c,\tilde{c}>0$ independent
of $\ell$ such that, for $\delta_{\ell}<\ell$,
\begin{equation}
\left\Vert \nabla\left(  u_{\ell}-\psi_{\ell}\left(  \lambda_{\infty}%
,\cdot\right)  \otimes u_{\infty}\right)  \right\Vert _{2,\Omega_{\ell
-\delta_{\ell}}}^{2}\leq C_{\omega,\delta_{\ell}}^{(1)}\left\Vert u_{\infty
}\right\Vert _{2,\omega}^{2}+C_{\omega,\delta_{\ell}}^{(2)}\left\Vert
\nabla^{\prime}u_{\infty}\right\Vert _{2,\omega}^{2}, \label{Error_psi}%
\end{equation}
with
\[
C_{\omega,\delta_{\ell}}^{(1)}:=4\operatorname{e}^{-2\lambda_{1}\delta_{\ell}%
}\text{,}\qquad C_{\omega,\delta_{\ell}}^{(2)}:=4\left(  \frac{1}{\lambda_{1}%
}\operatorname{e}^{-2\lambda_{1}\delta_{\ell}}+6\left(  \ell-\delta_{\ell
}\right)  \operatorname{e}^{-2\lambda_{1}\ell}\right)
\]
and%
\begin{equation}
\lambda_{1}:=\inf_{v\in H_{0}^{1}\left(  \omega\right)  \backslash\left\{
0\right\}  }\frac{\left\Vert \nabla^{\prime}v\right\Vert _{2,\omega}%
}{\left\Vert v\right\Vert _{2,\omega}}. \label{deflambda1}%
\end{equation}
The right-hand side in (\ref{Error_psi}) goes to 0 with an exponential rate of
convergence if $\delta_{\ell}$ is bounded from below when $\ell\rightarrow
\infty$.
\end{lemma}

%

\proof
For $i=1,2,\ldots$, let $w_{i}$ be the $i$-th eigenfunction of $-\Delta
^{\prime}$, i.e., $w_{i}\in H_{0}^{1}\left(  \omega\right)  $ is a solution
of
\begin{equation}
\left(  \nabla^{\prime}w_{i},\nabla^{\prime}v\right)  _{2,\omega}=\lambda
_{i}^{2}\left(  w_{i},v\right)  _{2,\omega}\qquad\forall v\in H_{0}^{1}\left(
\omega\right)  \label{defwi}%
\end{equation}
and we normalize the eigenfunctions such that $\left(  w_{i},w_{j}\right)
_{2,\omega}=\delta_{i,j}$ and order them such that $\left(  \lambda
_{i}\right)  _{i}$ is increasing monotonously. Furthermore let $u_{\ell,i}\in
H_{0}^{1}(\Omega_{\ell})$ be the solution of%
\[
\left(  \nabla u_{\ell,i},\nabla v\right)  _{2}=\lambda_{i}^{2}\left(
1\otimes w_{i},v\right)  _{2}\qquad\forall v\in H_{0}^{1}\left(  \Omega_{\ell
}\right)  .
\]
Then one concludes from (\ref{dglpsi}) and (\ref{defwi}) that
\begin{equation}
u_{\ell,i}=\psi_{\ell}\left(  \lambda_{i},\cdot\right)  \otimes w_{i}.
\label{DefAlpha}%
\end{equation}
If $f\in L^{2}(\omega)$ it holds
\[
f=\sum_{i=1}^{\infty}(f,w_{i})_{2,\omega}w_{i}.
\]
This shows that the solutions of (\ref{reducedDBVP}) and (\ref{originalDBVP})
can be expressed as
\begin{align*}
u_{\infty}  &  =\sum_{i=1}^{\infty}\frac{(f,w_{i})_{2,\omega}}{\lambda_{i}%
^{2}}w_{i},\\
u_{\ell}  &  =\sum_{i=1}^{\infty}\frac{(f,w_{i})_{2,\omega}}{\lambda_{i}^{2}%
}u_{\ell,i}=\sum_{i=1}^{\infty}\frac{(f,w_{i})_{2,\omega}}{\lambda_{i}}%
\psi_{\ell}\left(  \lambda_{i},\cdot\right)  \otimes w_{i}.
\end{align*}
With $\psi_{\ell}$ as in (\ref{approximationM1}) we get
\[
u_{\ell}-\psi_{\ell}\left(  \lambda_{\infty},\cdot\right)  \otimes u_{\infty
}=\sum_{i=1}^{\infty}\frac{(f,w_{i})_{2,\omega}}{\lambda_{i}^{2}}\phi_{\ell,i}\otimes
w_{i},\text{\qquad}\phi_{\ell,i}(x_{1}):=\frac{\cosh(\lambda_{\infty}x_{1}%
)}{\cosh(\lambda_{\infty}\ell)}-\frac{\cosh(\lambda_{i}x_{1})}{\cosh
(\lambda_{i}\ell)}.
\]
Let $\delta_{\ell}<\ell$. Then, since $\int_{\omega}w_{i}w_{j}dx^{\prime
}=\delta_{i,j}$, we get
\begin{align}
|\nabla\left(  u_{\ell}-\psi_{\ell}\left(  \lambda_{\infty},\cdot\right)
\otimes u_{\infty}\right)  |_{2,\Omega_{\ell-\delta_{\ell}}}^{2}  &
=\int_{-\ell+\delta_{\ell}}^{\ell-\delta_{\ell}}\int_{\omega}\left\vert
\nabla(u_{\ell}-\psi_{\ell}\left(  \lambda_{\infty},\cdot\right)  \otimes
u_{\infty})\right\vert ^{2}dx\nonumber\\
&  =\sum_{i=1}^{\infty}\frac{(f,w_{i})^{2}}{\lambda_{i}^{4}}\int_{-\ell
+\delta_{\ell}}^{\ell-\delta_{\ell}}\left(  \left(  \phi_{\ell,i}^{\prime
}\right)  ^{2}+\lambda_{i}^{2}\phi_{\ell,i}^{2}\right)  . \label{uluinfest}%
\end{align}
One has for any $\alpha>0$
\begin{align*}
\int_{-\ell+\delta_{\ell}}^{\ell-\delta_{\ell}}\left(  \frac{\cosh(\alpha
x_{1})}{\cosh(\alpha\ell)}\right)  ^{2}dx_{1}  &  =2\int_{0}^{\ell
-\delta_{\ell}}\left(  \frac{\cosh(\alpha x_{1})}{\cosh(\alpha\ell)}\right)
^{2}dx_{1}\\
&  =\frac{1}{2}\int_{0}^{\ell-\delta_{\ell}}\frac{\operatorname{e}^{2\alpha
x_{1}}+2+\operatorname{e}^{-2\alpha x_{1}}}{\cosh(\alpha\ell)^{2}}dx_{1}\\
&  \leq2\int_{0}^{\ell-\delta_{\ell}}\frac{\operatorname{e}^{2\alpha x_{1}}%
+3}{\operatorname{e}^{2\alpha\ell}}dx_{1}\\
&  \leq\frac{1}{\alpha}\operatorname{e}^{-2\alpha\delta_{\ell}}+6(\ell
-\delta_{\ell})\operatorname{e}^{-2\alpha\ell}%
\end{align*}
and similarly
\[
\int_{-\ell+\delta_{\ell}}^{\ell-\delta_{\ell}}\left(  \frac{\sinh(\alpha
x_{1})}{\cosh(\alpha\ell)}\right)  ^{2}dx_{1}\leq\frac{1}{\alpha
}\operatorname{e}^{-2\alpha\delta_{\ell}}.
\]
Since $\lambda_{1}\leq\lambda_{i}$ for all $i\in\mathbb{N}$ and
\[
\lambda_{1}=\inf_{v\in H_{0}^{1}\left(  \omega\right)  \backslash\left\{
0\right\}  }\frac{\left\Vert \nabla^{\prime}v\right\Vert _{2,\omega}%
}{\left\Vert v\right\Vert _{2,\omega}}\leq\frac{\left\Vert \nabla^{\prime
}u_{\infty}\right\Vert _{2,\omega}}{\left\Vert u_{\infty}\right\Vert
_{2,\omega}}=\lambda_{\infty}%
\]
we get
\begin{equation}
\int_{-\ell+\delta_{\ell}}^{\ell-\delta_{\ell}}\left(  \phi_{\ell,i}^{\prime
}\right)  ^{2}\leq2\operatorname{e}^{-2\lambda_{\infty}\delta_{\ell}%
}+2\operatorname{e}^{-2\lambda_{1}\delta_{\ell}}\leq4\operatorname{e}%
^{-2\lambda_{1}\delta_{\ell}}=C_{\omega,\delta_{\ell}}^{(1)}
\label{estphiprime}%
\end{equation}
and%
\begin{align}
\int_{-\ell+\delta_{\ell}}^{\ell-\delta_{\ell}}\phi_{\ell,i}^{2}  &  \leq
2\int_{-\ell+\delta_{\ell}}^{\ell-\delta_{\ell}}\left(  \left\vert \frac
{\cosh(\lambda_{\infty}x_{1})}{\cosh(\lambda_{\infty}\ell)}\right\vert
^{2}+\left\vert \frac{\cosh(\lambda_{i}x_{1})}{\cosh(\lambda_{i}\ell
)}\right\vert ^{2}\right) \nonumber\\
&  \leq2\left(  \frac{1}{\lambda_{\infty}}\operatorname{e}^{-2\lambda_{\infty
}\delta_{\ell}}+\frac{1}{\lambda_{i}}\operatorname{e}^{-2\lambda_{i}%
\delta_{\ell}}+6\left(  \ell-\delta_{\ell}\right)  \operatorname{e}%
^{-2\lambda_{\infty}\ell}+6\left(  \ell-\delta_{\ell}\right)  \operatorname{e}%
^{-2\lambda_{i}\ell}\right) \nonumber\\
&  \leq4\left(  \frac{1}{\lambda_{1}}\operatorname{e}^{-2\lambda_{1}%
\delta_{\ell}}+6\left(  \ell-\delta_{\ell}\right)  \operatorname{e}%
^{-2\lambda_{1}\ell}\right)  =C_{\omega,\delta_{\ell}}^{(2)}. \label{estphi}%
\end{align}
We employ the estimates (\ref{estphiprime}) and (\ref{estphi}) in
(\ref{uluinfest}) and obtain%
\begin{align*}
\left\Vert \nabla\left(  u_{\ell}-\psi_{\ell}\left(  \lambda_{\infty}%
,\cdot\right)  \otimes u_{\infty}\right)  \right\Vert _{2,\Omega_{\ell
-\delta_{\ell}}}^{2}  &  \leq C_{\omega,\delta_{\ell}}^{(1)}\sum_{i=1}%
^{\infty}\frac{(f,w_{i})_{2,\omega}^{2}}{\lambda_{i}^{4}}+C_{\omega,\delta_{\ell}}%
^{(2)}\sum_{i=1}^{\infty}\frac{(f,w_{i})_{2,\omega}^{2}}{\lambda_{i}^{2}}\\
&  =C_{\omega,\delta_{\ell}}^{(1)}\left\Vert u_{\infty}\right\Vert _{2,\omega
}^{2}+C_{\omega,\delta_{\ell}}^{(2)}\left\Vert \nabla^{\prime}u_{\infty
}\right\Vert _{2,\omega}^{2},
\end{align*}
which shows the assertion.%
\endproof

Lemma \ref{convergenceResultM1} suggests that one cannot expect convergence of
the approximation $\psi_{\ell}\left(  \lambda_{\infty},\cdot\right)  \otimes
u_{\infty}$ on the whole domain $\Omega_{\ell}$. Indeed it can be shown that,
in general, $\left\Vert \nabla\left(  u_{\ell}-\psi_{\ell}\left(
\lambda_{\infty},\cdot\right)  \otimes u_{\infty}\right)  \right\Vert
_{2,\Omega_{\ell}}\nrightarrow0$ as $\ell\rightarrow\infty$. Setting
$\delta_{\ell}=0$ in Lemma \ref{convergenceResultM1} shows that the error on
$\Omega_{\ell}$ can be estimated as follows:

\begin{corollary}
It holds%
\[
\left\Vert \nabla\left(  u_{\ell}-\psi_{\ell}\left(  \lambda_{\infty}%
,\cdot\right)  \otimes u_{\infty}\right)  \right\Vert _{2,\Omega_{\ell}}%
^{2}\leq4\left(  \left\Vert u_{\infty}\right\Vert _{2,\omega}^{2}%
+6\ell\operatorname{e}^{-2\lambda_{1}\ell}\left\Vert \nabla^{\prime}u_{\infty
}\right\Vert _{2,\omega}^{2}\right)  ,
\]
where $\lambda_{1}$ is as in (\ref{deflambda1}).
\end{corollary}

\subsection{Method 2: An alternating least squares type iteration}

Method 1 can be interpreted as a 2-step algorithm to obtain an approximation
$u^{M_{1}}_{\ell}$ of $u_{\ell}$.

\begin{itemize}
\item {Step 1:} Solve (\ref{reducedDBVP}) in order to obtain an approximation
of the form $1\otimes u_{\infty}$ which is \textit{non-conforming}, i.e., does
not belong to $H_{0}^{1}\left(  \Omega_{\ell}\right)  $.

\item {Step 2:} Using $u_{\infty}$, find a function $\psi_{\ell}$ that
satisfies (\ref{bestApproximation1}) in order to obtain the \textit{conforming
}approximation $u_{\ell}^{M_{1}}:=\psi_{\ell}\left(  \lambda_{\infty}%
,\cdot\right)  \otimes u_{\infty}\in H_{0}^{1}\left(  \Omega_{\ell}\right)  $.
\end{itemize}

In this section we extend this idea and seek approximations of the form
\begin{equation}
u_{\ell,m}^{M_{2}}=\sum_{j=0}^{m}p^{(j)}\otimes q^{(j)} \label{repuM2}%
\end{equation}
by iteratively solving least squares problems similar to
(\ref{bestApproximation1}). We denote by
\[
\operatorname*{Res}\nolimits_{m}=u_{\ell}-u_{\ell,m}^{M_{2}}=u_{\ell}%
-\sum_{j=0}^{m}p^{(j)}\otimes q^{(j)}%
\]
the residual of the approximation and suggest the following iteration to
obtain $u_{\ell,m}^{M_{2}}$:

\begin{itemize}

\item $m=0$: Set $q^{(0)}=u_{\infty}$ and $p^{(0)}=\psi_{\ell}\left(
\lambda_{\infty},\cdot\right)  $.

\item $m>0$: Find $q^{(m)}\in H_{0}^{1}(\omega)$ s.t.
\begin{equation}
q^{(m)}=\argmin_{q\in H_{0}^{1}(\omega)}\left\Vert \nabla\left(
\operatorname*{Res}\nolimits_{m-1}-p^{(m-1)}\otimes q\right)  \right\Vert
_{2}. \label{minimizeForQ}%
\end{equation}
Then, given $q^{(m)}$, find $p^{(m)}\in H_{0}^{1}(I_{\ell})$ s.t.%
\begin{equation}
p^{(m)}=\argmin_{p\in H_{0}^{1}(I_{\ell})}\left\Vert \nabla\left(
\operatorname*{Res}\nolimits_{m-1}-p\otimes q^{(m)}\right)  \right\Vert _{2}.
\label{minimizeForP}%
\end{equation}
Iterate \eqref{minimizeForQ} and \eqref{minimizeForP} until a stopping
criterion is reached (inner iteration). Then set $\operatorname*{Res}%
\nolimits_{m} = \operatorname*{Res}\nolimits_{m-1} - p^{(m)}\otimes q^{(m)}$.
\end{itemize}

\noindent The algorithm exhibits properties of a greedy algorithm. It is easy
to see that in each step of the (outer) iteration the error decreases or stays
constant. We focus here on its accuracy in comparison with the two other
methods via numerical experiments. We emphasize that for tensors of order at
least $3$, convergence can be shown for the (inner)
ALS iteration (see \cite{Uschmajew_15,espig2015convergence,oseledets2018alternating,
uschmajew2012local}). This limit, however, is not a global minimum in general. The outer iteration can be shown to converge as well against the true solution $u_\ell$ under the condition that we find the best rank-1 approximation in the inner iteration (see \cite{falco2012proper}).\medskip

\noindent The idea of computing approximations in the separated form \eqref{repuM2} by iteratively enriching the current solution with rank-1 terms is known in the literature as \emph{Proper Generalized Decomposition} (PGD). The PGD has been applied to various problems in computational mechanics (e.g. \cite{mokdad2007simulation, ammar2007new, aghighi2013non, chinesta2007reduction}), computational rheology (\cite{chinesta2011overview}), quantum chemistry (e.g. \cite{ammar2008nanometric, ammar2012reduction}) and others.\\
An extensive review of the method can be found in \cite{chinesta2013proper}. For an error and convergence analysis of the (outer) iteration in the case of the Poisson equation we also refer to \cite{le2009results}, where a similar (but not identical) approach as ours is considered.\medskip

\noindent In each step of the (outer) iteration above we need to solve at least two
minimization problems (\ref{minimizeForQ}) and (\ref{minimizeForP}). In the
following we derive the strong formulations of these problems.

\subsubsection{Resolution of (\ref{minimizeForQ})}

As before an investigation of the functional
\[
J(q^{(m)}):=\left\Vert \nabla(\operatorname*{Res}\nolimits_{m-1}%
-p^{(m-1)}\otimes q^{(m)})\right\Vert _{2}^{2}%
\]
shows that $q^{(m)}$ needs to satisfy%
\begin{align*}
&  \left(  \nabla\left(  p^{\left(  m-1\right)  }\otimes q^{\left(  m\right)
}\right)  ,\nabla\left(  p^{\left(  m-1\right)  }\otimes q\right)  \right)
_{2}=\left(  \nabla\operatorname*{Res}\nolimits_{m-1},\nabla\left(  p^{\left(
m-1\right)  }\otimes q\right)  \right)  _{2}\\
\iff &  p_{0,m-1}\left(  -\Delta^{\prime}q^{\left(  m\right)  },q\right)
_{2,\omega}+p_{1,m-1}\left(  q^{\left(  m\right)  },q\right)  _{2,\omega
}=\left(  -\Delta\operatorname*{Res}\nolimits_{m-1},p^{\left(  m-1\right)
}\otimes q\right)  _{2}%
\end{align*}
for all $q\in H_{0}^{1}(\omega)$, where
\[
p_{0,m-1}:=\Vert p^{(m-1)}\Vert_{2,I_{\ell}}^{2},\quad p_{1,m-1}:=\Vert\left(
p^{(m-1)}\right)  ^{\prime}\Vert_{2,I_{\ell}}^{2}.
\]
For the right-hand side we obtain
\begin{align*}
&  \left(  -\Delta\operatorname*{Res}\nolimits_{m-1},p^{(m-1)}\otimes
q\right)  _{2}=\left(  -\Delta\left(  u_{\ell}-\sum_{j=0}^{m-1}p^{(j)}\otimes
q^{(j)}\right)  ,p^{(m-1)}\otimes q\right)  _{2}\\
&  \qquad=\left(  1\otimes f,p^{(m-1)}\otimes q\right)  _{2}+\left(
\sum_{j=0}^{m-1}\left(  p^{(j)}\right)  ^{\prime\prime}\otimes q^{(j)}%
+p^{(j)}\otimes\Delta^{\prime}q^{(j)},p^{(m-1)}\otimes q\right)  _{2}\\
&  \qquad=\tilde{p}_{m-1}\left(  f,q\right)  _{2,\omega}+\sum_{j=0}%
^{m-1}\left(  \tilde{p}_{2,j,m-1}\left(  q^{(j)},q\right)  _{2,\omega}%
+\tilde{p}_{0,j,m-1}\left(  \Delta^{\prime}q^{(j)},q\right)  _{2,\omega
}\right)  ,
\end{align*}
where
\[
\tilde{p}_{m-1}:=\int_{-\ell}^{\ell}p^{(m-1)},\text{\ }\tilde{p}%
_{2,j,m-1}:=\left(  \left(  p^{(j)}\right)  ^{\prime\prime},p^{(m-1)}\right)
_{2,I_{\ell}},\text{\ }\tilde{p}_{0,j,m-1}:=\left(  p^{(j)},p^{(m-1)}\right)
_{2,I_{\ell}}.
\]
In order to compute (\ref{minimizeForQ}) we therefore have to solve in
$\omega$%
\begin{equation}
-p_{0,m-1}\Delta^{\prime}q^{(m)}+p_{1,m-1}q^{(m)}=\tilde{p}_{m-1}f+\sum
_{j=0}^{m-1}\left(  \tilde{p}_{2,j,m-1}q^{(j)}+\tilde{p}_{0,j,m-1}%
\Delta^{\prime}q^{(j)}\right)  .\label{strongFormQ}%
\end{equation}

\subsubsection{Resolution of (\ref{minimizeForP})}

Setting the derivative of the functional
\[
J(p^{(m)}):=\Vert\nabla(\operatorname*{Res}\nolimits_{m-1}-p^{(m)}\otimes
q^{(m)})\Vert_{2}^{2}%
\]
to zero, shows that $p^{(m)}$ needs to satisfy%
\begin{align*}
\left(  \nabla\left(  p^{\left(  m\right)  }\otimes q^{\left(  m\right)
}\right)  ,\nabla\left(  p\otimes q^{\left(  m\right)  }\right)  \right)
_{2}  &  =\left(  \nabla\operatorname*{Res}\nolimits_{m},\nabla\left(
p\otimes q^{\left(  m\right)  }\right)  \right)  _{2}\\
-q_{0,m}\left(  \left(  p^{\left(  m\right)  }\right)  ^{\prime\prime
},p\right)  _{2,I_{\ell}}+q_{1,m}\left(  p^{\left(  m\right)  },p\right)
_{2,I_{\ell}}  &  =\left(  -\Delta\operatorname*{Res}\nolimits_{m},p\otimes
q^{\left(  m\right)  }\right)  _{2}%
\end{align*}
for all $p\in H_{0}^{1}(-\ell,\ell)$, where
\[
q_{0,m}=\Vert q^{(m)}\Vert_{2,\omega}^{2},\quad q_{1,m}=\Vert\nabla^{\prime
}q^{(m)}\Vert_{2,\omega}^{2}.
\]
For the right-hand side we obtain
\begin{align*}
&  \left(  -\Delta\operatorname*{Res}\nolimits_{m-1},p\otimes q^{\left(
m\right)  }\right)  _{2}=\left(  -\Delta\left(  u_{\ell}-\sum_{j=0}%
^{m-1}p^{(j)}\otimes q^{(j)}\right)  ,p\otimes q^{\left(  m\right)  }\right)
_{2}\\
&  \qquad=\left(  1\otimes f,p\otimes q^{\left(  m\right)  }\right)
_{2}+\left(  \sum_{j=0}^{m-1}\left(  \left(  p^{(j)}\right)  ^{\prime\prime
}\otimes q^{(j)}+p^{(j)}\otimes\Delta^{\prime}q^{(j)}\right)  ,p\otimes
q^{\left(  m\right)  }\right)  _{2}\\
&  \qquad=\tilde{q}_{m}\int_{-\ell}^{\ell}p+\sum_{j=0}^{m-1}\left(  \tilde
{q}_{0,j,m}\left(  \left(  p^{(j)}\right)  ^{\prime\prime},p\right)
_{2,I_{\ell}}+\tilde{q}_{2,j,m}\left(  p^{(j)},p\right)  _{2,I_{\ell}}\right)
,
\end{align*}
where
\[
\tilde{q}_{m}:=\left(  f,q^{(m)}\right)  _{2,\omega},\quad\tilde{q}%
_{2,j,m}:=\left(  \Delta^{\prime}q^{(j)},q^{(m)}\right)  _{2,\omega}%
,\quad\tilde{q}_{0,j,m}:=\left(  q^{(j)},q^{(m)}\right)  _{2,\omega}.
\]
In order to obtain the solution of (\ref{minimizeForP}) we therefore have to
solve in $I_{\ell}$%
\begin{equation}
-q_{0,m}\left(  p^{(m)}\right)  ^{\prime\prime}+q_{1,m}p^{(m)}=\tilde{q}%
_{m}+\sum_{j=0}^{m-1}\left(  \tilde{q}_{2,j,m}p^{(j)}+\tilde{q}_{0,j,m}\left(
p^{(j)}\right)  ^{\prime\prime}\right)  . \label{strongFormP}%
\end{equation}

\begin{remark}
The constants $p_{1,m-1}$, $\tilde{p}_{2,j,m-1}$, $q_{1,m}$ and $\tilde
{q}_{2,j,m}$ involve derivatives and Laplace-operators. Note that after
solving (\ref{strongFormQ}) and (\ref{strongFormP}) for $q^{(m)}$ and
$p^{(m)}$, discrete versions of $\Delta^{\prime}q^{(m)}$ and $\left(
p^{(m)}\right)  ^{\prime\prime}$ can be easily obtained via the same
equations. Furthermore, since
\[
q_{1,m}=\Vert\nabla^{\prime}q^{(m)}\Vert_{2,\omega}^{2}=\left(  -\Delta
^{\prime}q^{(m)},q^{(m)}\right)  _{2,\omega}=-\tilde{q}_{2,m,m}%
\]
and
\[
p_{1,m-1}=\Vert\left(  p^{(m-1)}\right)  ^{\prime}\Vert_{2,I_{\ell}}%
^{2}=\left(  -\left(  p^{(m-1)}\right)  ^{\prime\prime},p^{(m-1)}\right)
_{2,I_{\ell}}=-\tilde{p}_{2,m-1,m-1}%
\]
a numerical computation of the gradients can be avoided.
\end{remark}

\subsection{Method 3: Exploiting the tensor product structure of the operator}

In this section we exploit the tensor product structure of the Laplace
operator and the domain $\Omega_{\ell}$. Recall that
\[
\Omega_{\ell}=I_{\ell}\times\omega.
\]
Note that we do not assume that $\omega$ has a tensor product structure.
Furthermore the Laplace operator in our original problem (\ref{originalDBVP})
can be written as
\begin{equation}
-\Delta=-\partial_{1}^{2}-\Delta^{\prime}. \label{defmLap}%
\end{equation}
We discretise (\ref{originalDBVP}) with $F$ as in (\ref{defFtensor1}) on a
mesh $\mathcal{G}$, e.g., by finite elements or finite differences on a tensor
mesh, i.e., each mesh cell has the form $\left(  x_{i-1},x_{i}\right)
\times\tau_{j}$, where $\tau_{j}$ is an element of the mesh for $\omega$. The
essential assumption is that the system matrix for the discrete version of
$-\Delta$ in (\ref{defmLap}) is of the tensor form
\begin{equation}
A=A_{1}\otimes M^{\prime}+M_{1}\otimes A^{\prime}. \label{DiscretizedLaplace}%
\end{equation}
If we discretise with a finite difference scheme on an equidistant grid for
$I_{\ell}$ with step size $h$, then $A_{1}$ is the tridiagonal matrix
$h^{-2}\operatorname*{tridiag}\left[  -1,2,-1\right]  $ and $M_{1}$ is the
identity matrix. A finite element discretisation with piecewise linear
elements leads as well to $A_{1}=h^{-2}\operatorname*{tridiag}\left[
-1,2,-1\right]  $, while $M_{x_{1}}=\operatorname*{tridiag}\left[
1/6,2/3,1/6\right]  $. It can be shown that the inverse of the matrix $A$ can
be efficiently approximated with a sum of matrix exponentials. More precisely
the following Theorem holds which is proved in \cite{hackbusch2012tensor},
Proposition 9.34.

\begin{theorem}
\label{ThmHackbusch}Let $M^{(j)}$, $A^{(j)}$ be positive definite matrices
with $\lambda_{\min}^{(j)}$ and $\lambda_{\max}^{(j)}$ being the extreme
eigenvalues of the generalized eigenvalue problem $A^{(j)}x=\lambda M^{(j)}x$
and set
\begin{align*}
A=  &  A^{(1)}\otimes M^{(2)}\otimes\ldots\otimes M^{(n)}+M^{(1)}\otimes
A^{(2)}\otimes\ldots\otimes M^{(n)}+\ldots\\
&  +M^{(1)}\otimes\ldots\otimes M^{(n-1)}\otimes A^{(n)}.
\end{align*}
Then $A^{-1}$ can be approximated by
\[
B:=\left(  \sum_{\nu=1}^{r}a_{\nu,[a,b]}\bigotimes_{j=1}^{n}\operatorname{exp}%
\left(  -\alpha_{\nu,[a,b]}\left(  M^{(j)}\right)  ^{-1}A^{(j)}\right)
\right)  \left(  \bigotimes_{j=1}^{n}\left(  M^{(j)}\right)  ^{-1}\right)  ,
\]
where the coefficients $a_{\nu},\alpha_{\nu}>0$ are such that
\begin{align*}
\varepsilon (\frac{1}{x},[a,b],r)  &  :=\left\Vert \frac{1}{x}-\sum_{\nu=1}%
^{r}a_{\nu,[a,b]}\operatorname{e}^{-\alpha_{\nu,[a,b]}x}\right\Vert _{[a,b]}\\
&  =\inf\left\{  \left\Vert \frac{1}{x}-\sum_{\nu=1}^{r}b_{\nu}%
\operatorname{e}^{-\beta_{\nu}x}\right\Vert _{[a,b],\infty}:b_{\nu},\beta
_{\nu}\in\mathbb{R}\right\}
\end{align*}
with $a:=\sum_{j=1}^{n}\lambda_{\min}^{(j)}$ and $b:=\sum_{j=1}^{n}%
\lambda_{\max}^{(j)}$. The error can be estimated by
\[
\Vert A^{-1}-B\Vert_{2}\leq\varepsilon(\frac{1}{x},[a,b],r)\Vert M^{-1}%
\Vert_{2},
\]
where $M=\otimes_{j=1}^{n}M^{(j)}$.
\end{theorem}

Theorem \ref{ThmHackbusch} shows how the inverse of matrices of the form
(\ref{DiscretizedLaplace}) can be approximated by sums of matrix exponentials.
It is based on the approximability of the function $1/x$ by sums of
exponentials in the interval $[a,b]$. We refer to
\cite{hackbusch2012tensor,hackbuschCVS} for details how to choose $r$ and the
coefficients $a_{\nu,[a,b]}$, $\alpha_{\nu,[a,b]}$ in order to reach a given
error tolerance $\varepsilon(\frac{1}{x},[a,b],r)$. Note that the interval
$[a,b]$ where $1/x$ needs to be approximated depends on the matrices $A^{(j)}$
and $M^{(j)}$. Thus, if $A$ changes $a$ and $b$ need to be recomputed which in
turn has an influence on the optimal choice of the parameters $a_{\nu,[a,b]}$
and $\alpha_{\nu,[a,b]}$. \newline

Numerical methods based on Theorem \ref{ThmHackbusch} can only be efficient if
the occurring matrix exponential can be evaluated at low cost. In our setting
we will need to compute the matrices $\operatorname{exp}\left(  -\alpha
_{\nu,[a,b]}M_{1}^{-1}A_{1}\right)  $ and $\operatorname{exp}\left(
-\alpha_{\nu,[a,b]}\left(  M^{\prime}\right)  ^{-1}A^{\prime}\right)  $. The
evaluation of the first matrix will typically be simpler. In the case where a
finite difference scheme is employed and $A_{1}$ is a tridiagonal Toeplitz
matrix while $M_{1}$ is the identity, the matrix exponential can be computed
by diagonalizing $A_{1}$, e.g., $A_{1}=SD_{1}S^{-1}$, and using
$\operatorname{exp}\left(  -\alpha_{\nu,[a,b]}M_{1}^{-1}A_{1}\right)
=S\operatorname{exp}\left(  -\alpha_{\nu,[a,b]}D_{1}\right)  S^{-1}$. The
computation of exponentials for general matrices is more involved. We refer to
\cite{moler2003nineteen} for an overview of different numerical methods. Here,
we will make use of the Dunford-Cauchy integral (see
\cite{hackbusch2015hierarchical}). For a matrix $\tilde{M}$ we can write
\begin{equation}
\operatorname{exp}\left(  -\tilde{M}\right)  =\frac{1}{2\pi\operatorname{i}%
}\oint_{\mathcal{C}}\left(  \zeta I-\tilde{M}\right)  ^{-1}\operatorname{e}%
^{-\zeta}d\zeta\label{DunfordCauchy}%
\end{equation}
for a contour $\mathcal{C}=\partial D$ which encircles all eigenvalues of
$\tilde{M}$. We assume here that $\tilde{M}$ is positive definite. Then the
spectrum of $\tilde{M}$ satisfies $\sigma(\tilde{M})\subset(0,\Vert M\Vert]$
and the following (infinite) parabola
\[
\left\{  \zeta(s)=x(s)+\operatorname{i}y(s):x(s):=s^{2},y(s):=-s\quad\text{for
}s\in\mathbb{R}\right\}
\]
can be used as integration curve $\mathcal{C}$. The substitution
$\zeta\rightarrow s^{2}-\operatorname{i}s$ then leads to
\begin{equation}
\operatorname{exp}\left(  -\tilde{M}\right)  =\int_{-\infty}^{\infty
}\underbrace{\left(  \frac{1}{2\pi\operatorname{i}}(s^{2}-\operatorname{i}%
s)I-\tilde{M}\right)  ^{-1}\operatorname{e}^{-s^{2}+\operatorname{i}%
s}(2s-\operatorname{i})}_{=:G(s)}ds. \label{DunfordCauchySubs}%
\end{equation}
The integrand decays exponentially for $s\rightarrow\pm\infty$. Therefore
(\ref{DunfordCauchySubs}) can be efficiently approximated by sinc quadrature,
i.e.,
\begin{equation}
\operatorname{exp}\left(  -\tilde{M}\right)  =\int_{-\infty}^{\infty
}G(s)ds\approx\mathfrak{h}\sum_{\nu=-N}^{N}G(\nu\mathfrak{h}),
\label{DunfordCauchyApprox}%
\end{equation}
where $\mathfrak{h}>0$ and should be chosen s.t. $\mathfrak{h}=\mathcal{O}%
\left(  (N+1)^{-2/3}\right)  $. We refer to \cite{hackbusch2015hierarchical}
for an introduction to sinc quadrature and for error estimates for the
approximation in (\ref{DunfordCauchyApprox}). The parameters $\mathfrak{h}$
and $N$ in our implementation have been chosen such that quadrature errors
become negligible compared to the overall discretisation error. For practical
computations, the \textit{halving rule} (see \cite[\S 14.2.2.2]%
{hackbusch2015hierarchical}) could be faster while the Dunford-Schwartz
representation with sinc quadrature is more suited for an error analysis.


\section{Numerical Experiments\label{NumericalExperiments}}

\subsection{The case of a planar cylinder}

In this subsection we apply the methods derived in Section
\ref{NumericlaApproximation} to a simple model problem in two dimensions. We
consider the planar cylinder
\[
\Omega_{\ell}^{2D}=I_{\ell}\times(-1,1)
\]
and solve (\ref{originalDBVP}) for different right-hand sides $F=1\otimes f$
(see (\ref{defFtensor1})) and different lengths $\ell$. The reduced problem
(\ref{reducedDBVP}) on $\omega=(-1,1)$ is solved using a standard finite
difference scheme. We compare the approximations of (\ref{originalDBVP}) to a
reference solution $u_{2D,\ell}^{\operatorname*{ref}}$ that is computed using
a finite difference method on sufficiently refined two-dimensional
grid.\newline In Table \ref{TableErrorsM1} we state the $L^{2}\left(
\Omega_{\ell}^{2D}\right)  $-errors of the approximations $u_{2D,\ell}^{M_{1}%
}$ for various values of $\ell$ and right-hand sides $f$. Having in mind that
$u_{2D,\ell}^{M_{1}}$ is a rather simple one-term approximation that only
requires the solution of one $(n-1)$-dimensional problem (plus some
postprocessing), the accuracy of the approximation is satisfactory especially
for larger values of $\ell$.

\begin{table}[ptb]
\centering
\begin{tabular}
[c]{lllll}
& $f(x^{\prime})=1$ & $f(x^{\prime})=\sin(2x^{\prime}+0.5)$ & $f(x^{\prime
})=\tanh(4x^{\prime}+1)$ & $f(x^{\prime})=|x^{\prime}|$\\\hline
\multicolumn{1}{|l|}{$\ell=1$} & \multicolumn{1}{c|}{$3.30\cdot10^{-2}$} &
\multicolumn{1}{c|}{$2.44\cdot10^{-1}$} & \multicolumn{1}{c|}{$2.48\cdot
10^{-1}$} & \multicolumn{1}{c|}{$1.15\cdot10^{-1}$}\\\hline
\multicolumn{1}{|l|}{$\ell=5$} & \multicolumn{1}{c|}{$6.34\cdot10^{-3}$} &
\multicolumn{1}{c|}{$5.24\cdot10^{-2}$} & \multicolumn{1}{c|}{$5.39\cdot
10^{-2}$} & \multicolumn{1}{c|}{$2.22\cdot10^{-2}$}\\\hline
\multicolumn{1}{|l|}{$\ell=10$} & \multicolumn{1}{c|}{$4.24\cdot10^{-3}$} &
\multicolumn{1}{c|}{$3.53\cdot10^{-2}$} & \multicolumn{1}{c|}{$3.63\cdot
10^{-2}$} & \multicolumn{1}{c|}{$1.51\cdot10^{-2}$}\\\hline
\multicolumn{1}{|l|}{$\ell=20$} & \multicolumn{1}{c|}{$2.92\cdot10^{-3}$} &
\multicolumn{1}{c|}{$2.44\cdot10^{-2}$} & \multicolumn{1}{c|}{$2.51\cdot
10^{-2}$} & \multicolumn{1}{c|}{$1.03\cdot10^{-2}$}\\\hline
\multicolumn{1}{|l|}{$\ell=50$} & \multicolumn{1}{c|}{$1.82\cdot10^{-3}$} &
\multicolumn{1}{c|}{$1.52\cdot10^{-2}$} & \multicolumn{1}{c|}{$1.56\cdot
10^{-2}$} & \multicolumn{1}{c|}{$6.45\cdot10^{-3}$}\\\hline
\end{tabular}
\caption{Relative $L^{2}\left(  \Omega_{\ell}^{2D}\right)  $-errors of the
approximations $u_{2D,\ell}^{M_{1}}$ for different values of $\ell$ and $f$.}%
\label{TableErrorsM1}%
\end{table}

\begin{figure}[tbh]
\centering
\includegraphics[width=0.95\textwidth]{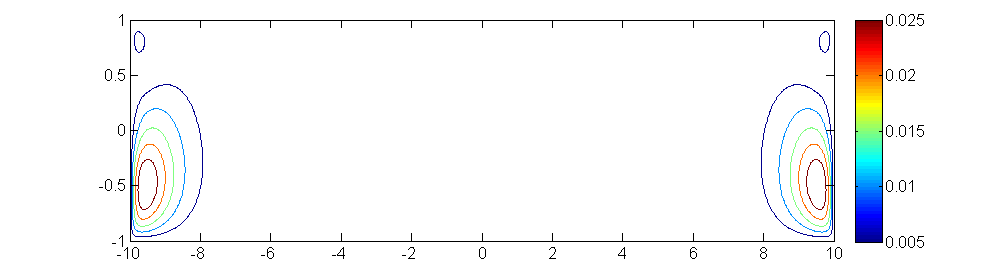}
\caption{Absolute error $|u_{2D,\ell}^{M_{1}}-u_{2D,\ell}^{\operatorname*{ref}%
}|$ for $\ell=10$ and $f(x^{\prime})=\tanh(4x^{\prime}+1)$.}%
\label{errorEll10}%
\end{figure}

Figure \ref{errorEll10} shows the pointwise, absolute error $|u_{2D,\ell
}^{M_{1}}-u_{2D,\ell}^{\operatorname*{ref}}|$ in $\Omega_{\ell}$ for $\ell=10$
and $f(x^{\prime})=\tanh(4x^{\prime}+1)$. As expected the accuracy of the
approximation is very high in the interior of the planar cylinder (away from
$\pm\ell$).

\begin{figure}[tbh]
\centering
\subfigure[$\ell = 20$]{
\centering
\includegraphics[width=0.47\textwidth]{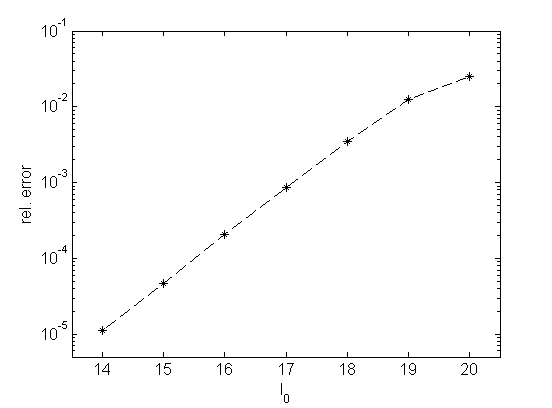}
} \hfil
\subfigure[$\ell = 50$ ]{
\centering
\includegraphics[width=0.47\textwidth]{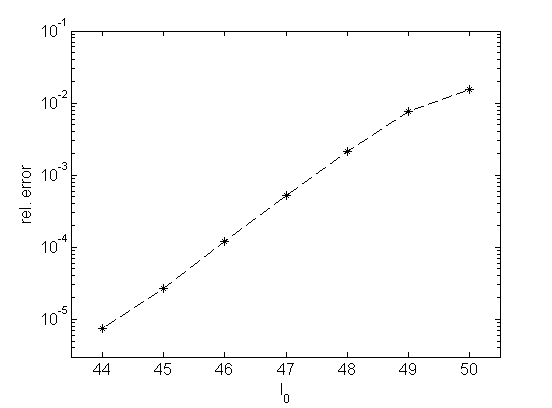}
} \caption{Relative $L^{2}$-errors of the approximation $u_{2D,\ell}^{M_{1}}$
in $\Omega_{\ell_{0}}$ for $f(x^{\prime})=\tanh(4x^{\prime}+1)$.}%
\label{RelativeErrors1D_wrt_l0}%
\end{figure}

Lemma \ref{convergenceResultM1} (and Figure \ref{errorEll10}) suggests that
the approximation in the interior of the cylinder is significantly better than
on the whole domain $\Omega_{\ell}$. Indeed, if the region of interest is only
a subdomain $\Omega_{\ell_{0}}\subset\Omega_{\ell}$, where $\ell_{0}<\ell$,
the error decreases exponentially as $\ell_{0}\rightarrow0$. Figure
\ref{RelativeErrors1D_wrt_l0} shows the relative error $\Vert u_{2D,\ell
}^{M_{1}}-u_{2D,\ell}^{\operatorname*{ref}}\Vert_{L^{2}(\Omega_{\ell_{0}}%
)}/\Vert u_{2D,\ell}^{\operatorname*{ref}}\Vert_{L^{2}(\Omega_{\ell_{0}})}$
with respect to $\ell_{0}$ for $\ell=20,50$ and the right-hand side
$f(x^{\prime})=\tanh(4x^{\prime}+1)$. We can see that the exponential
convergence sets in almost immediately as $l_{0}$ moves away from $\ell$.

\vspace{3mm} \noindent To conclude, Method 1 can be used in
applications where

\begin{itemize}
\item only a limited approximation accuracy is required,

\item a good starting point for more accurate methods is needed,

\item the region of interest is a subdomain $\Omega_{\ell_{0}}$ of
$\Omega_{\ell}$ with $\ell_{0}<\ell$.
\end{itemize}

\vspace{3mm} \noindent In Method 2 we use $u_{2D,\ell}^{M_{1}}$ as starting
value of the iteration which is then successively refined by approximating the
residual in each step with a series of $L^{2}$ best approximations. In Table
\ref{TableErrorsM2} we state the relative errors of this approach in the case
$f(x^{\prime})=\tanh(4x^{\prime}+1)$ for different values of $\ell$ and
iteration steps. We can see that five iterations are sufficient to reduce the
error of the initial approximation $u_{2D,\ell}^{M_{1}}$ by a factor $100$ for
all considered values of $\ell$. However, in this case more iterations do not
lead to significantly better results and the convergence seems to flatten. One
explanation for this is that the residuals are increasingly difficult to
approximate with each step of the iteration. After a few iterations a one-term
approximation of these residuals of the form $p^{(m)}\otimes q^{(m)}$
therefore is not sufficiently accurate which leads to reduced decay of the
error in the overall scheme. \newline Note that in the case $\ell=1$,
$\Omega_{\ell}$ cannot be considered as a \textquotedblleft
long\textquotedblright\ domain. Therefore, the initial approximation
$u_{2D,\ell}^{M_{1}}$ only exhibits a low accuracy. Nevertheless the error of
$u_{2D,\ell,m}^{M_{2}}$ decays quickly as $m$ increases and reaches a similar
level of accuracy as for larger $\ell$. This shows that Method 2 can be used for more general domains than considered here (e.g. \cite{giner2013proper}).

\begin{table}[ptb]
\centering
\begin{tabular}
[c]{cccccc}
& $\ell=1$ & $\ell=5$ & $\ell=10$ & $\ell=20$ & $\ell=50$\\\hline
\multicolumn{1}{|l|}{$m=1$} & \multicolumn{1}{c|}{$2.48\cdot10^{-1}$} &
\multicolumn{1}{c|}{$5.39\cdot10^{-2}$} & \multicolumn{1}{c|}{$3.63\cdot
10^{-2}$} & \multicolumn{1}{c|}{$2.51\cdot10^{-2}$} &
\multicolumn{1}{c|}{$1.57\cdot10^{-2}$}\\\hline
\multicolumn{1}{|l|}{$m=2$} & \multicolumn{1}{c|}{$1.96\cdot10^{-2}$} &
\multicolumn{1}{c|}{$1.32\cdot10^{-2}$} & \multicolumn{1}{c|}{$9.33\cdot
10^{-3}$} & \multicolumn{1}{c|}{$6.58\cdot10^{-3}$} &
\multicolumn{1}{c|}{$8.16\cdot10^{-3}$}\\\hline
\multicolumn{1}{|l|}{$m=3$} & \multicolumn{1}{c|}{$7.44\cdot10^{-3}$} &
\multicolumn{1}{c|}{$2.66\cdot10^{-3}$} & \multicolumn{1}{c|}{$1.85\cdot
10^{-3}$} & \multicolumn{1}{c|}{$1.29\cdot10^{-3}$} &
\multicolumn{1}{c|}{$8.16\cdot10^{-4}$}\\\hline
\multicolumn{1}{|l|}{$m=4$} & \multicolumn{1}{c|}{$1.18\cdot10^{-3}$} &
\multicolumn{1}{c|}{$7.73\cdot10^{-4}$} & \multicolumn{1}{c|}{$5.46\cdot
10^{-4}$} & \multicolumn{1}{c|}{$3.87\cdot10^{-4}$} &
\multicolumn{1}{c|}{$2.44\cdot10^{-4}$}\\\hline
\multicolumn{1}{|l|}{$m=5$} & \multicolumn{1}{c|}{$3.80\cdot10^{-4}$} &
\multicolumn{1}{c|}{$3.74\cdot10^{-4}$} & \multicolumn{1}{c|}{$2.71\cdot
10^{-4}$} & \multicolumn{1}{c|}{$1.96\cdot10^{-4}$} &
\multicolumn{1}{c|}{$1.25\cdot10^{-4}$}\\\hline
\multicolumn{1}{|l|}{$m=6$} & \multicolumn{1}{c|}{$1.68\cdot10^{-4}$} &
\multicolumn{1}{c|}{$3.05\cdot10^{-4}$} & \multicolumn{1}{c|}{$2.23\cdot
10^{-4}$} & \multicolumn{1}{c|}{$1.63\cdot10^{-4}$} &
\multicolumn{1}{c|}{$1.04\cdot10^{-4}$}\\\hline
\multicolumn{1}{|l|}{$m=7$} & \multicolumn{1}{c|}{$1.33\cdot10^{-4}$} &
\multicolumn{1}{c|}{$2.90\cdot10^{-4}$} & \multicolumn{1}{c|}{$2.12\cdot
10^{-4}$} & \multicolumn{1}{c|}{$1.55\cdot10^{-4}$} &
\multicolumn{1}{c|}{$9.95\cdot10^{-5}$}\\\hline
\end{tabular}
\caption{Relative $L^{2}$-errors of the approximations $u_{2D,\ell,m}^{M_{2}}$
for different values of $\ell$ and iterations $m$. We used $f(x^{\prime
})=\tanh(4x^{\prime}+1)$ throughout.}%
\label{TableErrorsM2}%
\end{table}

In Table \ref{TableErrorsM3} we show the relative errors of the approximations
$u_{2D,\ell,r}^{M_{3}}$ for $f(x^{\prime})=\tanh(4x^{\prime}+1)$ and different
values of $\ell$ and $r$. As the theory predicts the error decays
exponentially in $r$ and is governed by the approximability of the function
$1/x$ by exponential sums. Note that in this two-dimensional example the
arising matrix exponentials could be computed via diagonalization of the
involved finite difference matrices. An approximation of the Dunford-Cauchy
integral was not necessary in this case.

\begin{table}[ptb]
\centering
\begin{tabular}
[c]{ccccc}
& $\ell=1$ & $\ell=5$ & $\ell=10$ & $\ell=20$\\\hline
\multicolumn{1}{|l|}{$r=1$} & \multicolumn{1}{c|}{$1.48\cdot10^{-1}$} &
\multicolumn{1}{c|}{$1.02\cdot10^{-1}$} & \multicolumn{1}{c|}{$9.76\cdot
10^{-2}$} & \multicolumn{1}{c|}{$9.60\cdot10^{-2}$}\\\hline
\multicolumn{1}{|l|}{$r=2$} & \multicolumn{1}{c|}{$2.87\cdot10^{-2}$} &
\multicolumn{1}{c|}{$3.32\cdot10^{-2}$} & \multicolumn{1}{c|}{$3.14\cdot
10^{-2}$} & \multicolumn{1}{c|}{$3.09\cdot10^{-3}$}\\\hline
\multicolumn{1}{|l|}{$r=3$} & \multicolumn{1}{c|}{$9.06\cdot10^{-3}$} &
\multicolumn{1}{c|}{$9.78\cdot10^{-3}$} & \multicolumn{1}{c|}{$9.88\cdot
10^{-3}$} & \multicolumn{1}{c|}{$9.79\cdot10^{-3}$}\\\hline
\multicolumn{1}{|l|}{$r=4$} & \multicolumn{1}{c|}{$2.07\cdot10^{-3}$} &
\multicolumn{1}{c|}{$3.06\cdot10^{-3}$} & \multicolumn{1}{c|}{$2.81\cdot
10^{-3}$} & \multicolumn{1}{c|}{$2.70\cdot10^{-3}$}\\\hline
\multicolumn{1}{|l|}{$r=5$} & \multicolumn{1}{c|}{$1.24\cdot10^{-3}$} &
\multicolumn{1}{c|}{$1.11\cdot10^{-3}$} & \multicolumn{1}{c|}{$1.16\cdot
10^{-3}$} & \multicolumn{1}{c|}{$1.16\cdot10^{-3}$}\\\hline
\end{tabular}
\caption{Relative $L^{2}$-errors of the approximations $u_{2D,\ell,r}^{M_{3}}$
for different values of $\ell$ and $r$. We used $f(x^{\prime})=\tanh
(4x^{\prime}+1)$ throughout.}%
\label{TableErrorsM3}%
\end{table}

\subsection{A three-dimensional domain with a non-rectangular cross section}

\begin{figure}[tbh]
\centering
\subfigure[$\Omega_{\ell}$]{
\centering
\includegraphics[width=0.47\textwidth]{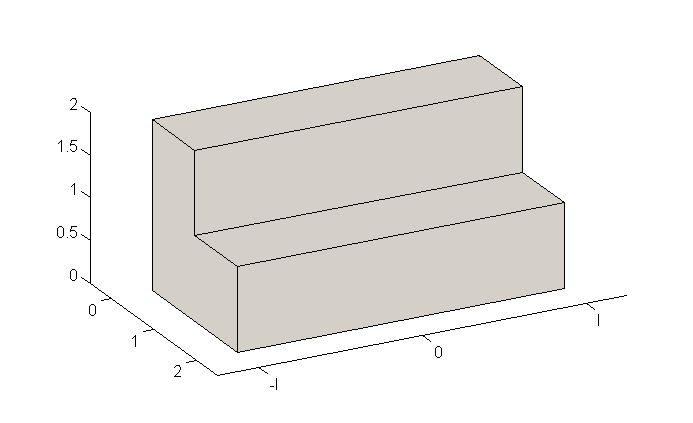}
} \hfil
\subfigure[$\omega$ ]{
\centering
\includegraphics[width=0.47\textwidth]{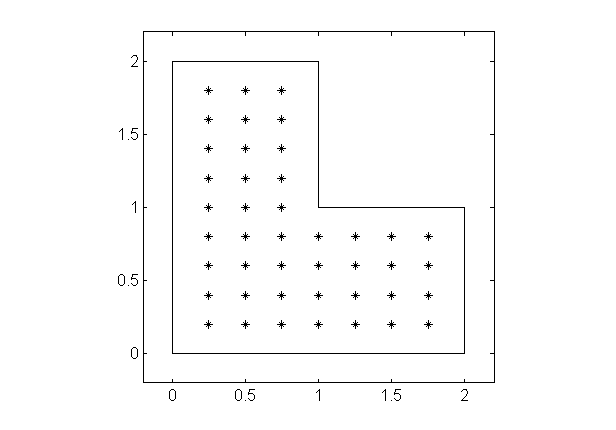}
} \caption{Plot of domain $\Omega_{\ell}$ and cross-section $\omega$.}%
\label{BoxShape}%
\end{figure}

In this section we consider the three-dimensional domain
\[
\Omega_{\ell}=(-\ell,\ell)\times\underbrace{\left[  (0,2)\times(0,1)\cup
(0,1)\times(1,2)\right]  }_{=\omega}%
\]
where $\omega$ is an ``L-shaped'' domain (see Figure \ref{BoxShape}). As before we
solve problem \eqref{originalDBVP} for different right-hand sides $f$ and
different values of $\ell$. The reduced problem \eqref{reducedDBVP} on
$\omega$ is solved using a standard 2D finite difference scheme. As 3D
reference solution we use an accurate approximation using method 3, i.e.
$u_{3D,\ell,r}^{M_{3}}$ for $r=30$, which is known to converge exponentially
in $r$.

\begin{table}[ptb]
\centering
\begin{tabular}
[c]{llll}
& $f(x^{\prime})=1$ & $f(x^{\prime})=\sin(x_{2}+0.5)x_{3}$ &
$f(x^{\prime})=\tanh(x_{2}x_{3})$\\\hline
\multicolumn{1}{|l|}{$\ell=1$} & \multicolumn{1}{c|}{$3.49\cdot10^{-2}$} &
\multicolumn{1}{c|}{$5.52\cdot10^{-2}$} & \multicolumn{1}{c|}{$5.59\cdot
10^{-2}$}\\\hline
\multicolumn{1}{|l|}{$\ell=5$} & \multicolumn{1}{c|}{$1.12\cdot10^{-2}$} &
\multicolumn{1}{c|}{$1.77\cdot10^{-2}$} & \multicolumn{1}{c|}{$1.79\cdot
10^{-2}$}\\\hline
\multicolumn{1}{|l|}{$\ell=10$} & \multicolumn{1}{c|}{$7.73\cdot10^{-3}$} &
\multicolumn{1}{c|}{$1.22\cdot10^{-2}$} & \multicolumn{1}{c|}{$1.23\cdot
10^{-2}$}\\\hline
\multicolumn{1}{|l|}{$\ell=20$} & \multicolumn{1}{c|}{$5.37\cdot10^{-3}$} &
\multicolumn{1}{c|}{$8.51\cdot10^{-3}$} & \multicolumn{1}{c|}{$8.86\cdot
10^{-2}$}\\\hline
\multicolumn{1}{|l|}{$\ell=50$} & \multicolumn{1}{c|}{$3.37\cdot10^{-3}$} &
\multicolumn{1}{c|}{$5.34\cdot10^{-3}$} & \multicolumn{1}{c|}{$5.54\cdot
10^{-2}$}\\\hline
\end{tabular}
\caption{Relative $L^{2}$-errors of the approximations $u_{3D,\ell}^{M_{1}}$
for different values of $\ell$ and $f$.}%
\label{TableErrorsM1_3D}%
\end{table}

Table \ref{TableErrorsM1_3D} shows the relative errors of the approximations
$u_{3D,\ell}^{M_{1}}$ for different values of $\ell$ and right-hand sides $f$.
As the theory predicts we cannot observe an exponentially decreasing error as
$\ell$ gets large, since we measure the error on the whole domain
$\Omega_{\ell}$ and not only a subdomain $\Omega_{\ell-\delta_{\ell}}$. As
before we only have to solve one two-dimensional problem on $\omega$ in order
to obtain the approximation $u_{3D,\ell}^{M_{1}}$.

\begin{table}[ptb]
\centering
\begin{tabular}
[c]{cccccc}
& $\ell=1$ & $\ell=5$ & $\ell=10$ & $\ell=20$ & $\ell=50$\\\hline
\multicolumn{1}{|l|}{$m=1$} & \multicolumn{1}{c|}{$5.59\cdot10^{-2}$} &
\multicolumn{1}{c|}{$1.79\cdot10^{-2}$} & \multicolumn{1}{c|}{$1.24\cdot
10^{-2}$} & \multicolumn{1}{c|}{$8.63\cdot10^{-3}$} &
\multicolumn{1}{c|}{$5.42\cdot10^{-3}$}\\\hline
\multicolumn{1}{|l|}{$m=2$} & \multicolumn{1}{c|}{$9.82\cdot10^{-3}$} &
\multicolumn{1}{c|}{$4.38\cdot10^{-3}$} & \multicolumn{1}{c|}{$3.08\cdot
10^{-3}$} & \multicolumn{1}{c|}{$2.18\cdot10^{-3}$} &
\multicolumn{1}{c|}{$1.37\cdot10^{-3}$}\\\hline
\multicolumn{1}{|l|}{$m=3$} & \multicolumn{1}{c|}{$2.86\cdot10^{-3}$} &
\multicolumn{1}{c|}{$1.09\cdot10^{-3}$} & \multicolumn{1}{c|}{$7.64\cdot
10^{-4}$} & \multicolumn{1}{c|}{$5.37\cdot10^{-4}$} &
\multicolumn{1}{c|}{$3.39\cdot10^{-4}$}\\\hline
\multicolumn{1}{|l|}{$m=4$} & \multicolumn{1}{c|}{$8.03\cdot10^{-4}$} &
\multicolumn{1}{c|}{$3.37\cdot10^{-4}$} & \multicolumn{1}{c|}{$2.37\cdot
10^{-4}$} & \multicolumn{1}{c|}{$1.63\cdot10^{-4}$} &
\multicolumn{1}{c|}{$1.03\cdot10^{-4}$}\\\hline
\multicolumn{1}{|l|}{$m=5$} & \multicolumn{1}{c|}{$3.46\cdot10^{-4}$} &
\multicolumn{1}{c|}{$1.38\cdot10^{-4}$} & \multicolumn{1}{c|}{$9.93\cdot
10^{-5}$} & \multicolumn{1}{c|}{$7.08\cdot10^{-5}$} &
\multicolumn{1}{c|}{$4.50\cdot10^{-5}$}\\\hline
\multicolumn{1}{|l|}{$m=6$} & \multicolumn{1}{c|}{$2.73\cdot10^{-4}$} &
\multicolumn{1}{c|}{$1.03\cdot10^{-4}$} & \multicolumn{1}{c|}{$7.50\cdot
10^{-5}$} & \multicolumn{1}{c|}{$5.37\cdot10^{-5}$} &
\multicolumn{1}{c|}{$3.42\cdot10^{-5}$}\\\hline
\multicolumn{1}{|l|}{$m=7$} & \multicolumn{1}{c|}{$2.59\cdot10^{-4}$} &
\multicolumn{1}{c|}{$9.51\cdot10^{-5}$} & \multicolumn{1}{c|}{$6.92\cdot
10^{-5}$} & \multicolumn{1}{c|}{$4.96\cdot10^{-5}$} &
\multicolumn{1}{c|}{$3.16\cdot10^{-5}$}\\\hline
\end{tabular}
\caption{Relative $L^{2}$-errors of the approximations $u_{3D,\ell,m}^{M_{2}}$
for different values of $\ell$ and iterations $m$. We used $f(x^{\prime
})=\tanh(x_{2}x_{3})$ throughout.}%
\label{TableErrorsM2_3D}%
\end{table}

In Table \ref{TableErrorsM2_3D} we show the relative errors of the
approximations $u_{3D,\ell,m}^{M_{2}}$ for $f(x^{\prime})=\tanh(x_{2}x_{3})$ and different values of $\ell$ and $m$ (number of
iterations). As in the 2D case this method significantly improves the initial
approximation $u_{3D,\ell,1}^{M_{2}}=u_{3D,\ell}^{M_{1}}$ using the
alternating least squares type iteration. However, also here we observe that
the convergence slows down when a certain accuracy is reached. We remark that
a good starting point for the iteration is crucial for this method. In all our
experiments $u_{3D,\ell}^{M_{1}}$ was a good choice which leads to a
convergence behaviour similar to the ones in Table \ref{TableErrorsM2_3D}.
Other choices often did not lead to satisfactory results.

\begin{table}[ptb]
\centering
\begin{tabular}
[c]{ccccc}
& $\ell=1$ & $\ell=5$ & $\ell=10$ & $\ell=20$\\\hline
\multicolumn{1}{|l|}{$r=1$} & \multicolumn{1}{c|}{$1.10\cdot10^{-1}$} &
\multicolumn{1}{c|}{$1.06\cdot10^{-1}$} & \multicolumn{1}{c|}{$1.06\cdot
10^{-1}$} & \multicolumn{1}{c|}{$1.06\cdot10^{-1}$}\\\hline
\multicolumn{1}{|l|}{$r=2$} & \multicolumn{1}{c|}{$2.18\cdot10^{-2}$} &
\multicolumn{1}{c|}{$1.94\cdot10^{-2}$} & \multicolumn{1}{c|}{$1.92\cdot
10^{-2}$} & \multicolumn{1}{c|}{$1.92\cdot10^{-2}$}\\\hline
\multicolumn{1}{|l|}{$r=3$} & \multicolumn{1}{c|}{$6.21\cdot10^{-3}$} &
\multicolumn{1}{c|}{$6.63\cdot10^{-3}$} & \multicolumn{1}{c|}{$6.63\cdot
10^{-3}$} & \multicolumn{1}{c|}{$6.62\cdot10^{-3}$}\\\hline
\multicolumn{1}{|l|}{$r=4$} & \multicolumn{1}{c|}{$2.36\cdot10^{-3}$} &
\multicolumn{1}{c|}{$2.02\cdot10^{-3}$} & \multicolumn{1}{c|}{$2.02\cdot
10^{-3}$} & \multicolumn{1}{c|}{$2.01\cdot10^{-3}$}\\\hline
\multicolumn{1}{|l|}{$r=5$} & \multicolumn{1}{c|}{$7.91\cdot10^{-4}$} &
\multicolumn{1}{c|}{$7.39\cdot10^{-4}$} & \multicolumn{1}{c|}{$7.15\cdot
10^{-4}$} & \multicolumn{1}{c|}{$7.15\cdot10^{-4}$}\\\hline
\end{tabular}
\caption{Relative $L^{2}$-errors of the approximations $u_{3D,\ell,r}^{M_{3}}$
for different values of $\ell$ and $r$. We used $f(x^{\prime})=\tanh
(x_{2}x_{3})$ throughout.}%
\label{TableErrorsM3_3D}%
\end{table}

In Table \ref{TableErrorsM3_3D} we show the relative errors of the
approximations $u_{3D,\ell,r}^{M_{3}}$ again for $f(x^{\prime})=\tanh
(x_{2}x_{3})$ and different values of $\ell$ and $r$. As
before the error decays exponentially with respect to $r$. The arising matrix
exponentials $\operatorname{exp}\left(  -\alpha_{\nu,[a,b]}A_{x^{\prime}%
}\right)  $ in these experiments were computed using the sinc quadrature
approximation \eqref{DunfordCauchyApprox}. The number of quadrature points $N$
was chosen such that the corresponding quadrature error had an negligible
effect on the overall approximation.

\section{Conclusion}

We have presented three different methods for constructing tensor
approximations to the solution of a Poisson equation on a long product domain
for a right-hand side which is an elementary tensor.

The construction of a one-term tensor approximation is based on asymptotic
analysis. The approximation converges exponentially (on a fixed subdomain) as
the length of the cylinder goes to infinity. However, the error is fixed for
fixed length since the approximation consists of only one term. The cost for
computing this approximation is very low -- it consists of solving a
Poisson-type problem on the cross section and a cheap post-processing step to
find the univariate function in the one-term tensor approximation.

The ALS type method uses this elementary tensor and generates step-by-step a
rank-$k$ approximation. The computation of the $m$--th term in the tensor
approximation itself requires an inner iteration. If one is interested in only
a moderate accuracy (but improved accuracy compared to the initial
approximation) this method is still relatively cheap and significantly
improves the accuracy. However, the theory for ALS for this application is not
fully developed and the definition of a good stopping criterion is based on
heuristics and experiments.

Finally the approximation which is based on exponential sums is the method of
choice if a higher accuracy is required. A well developed a priori error
analysis allows us to choose the tensor rank in the approximation in a very
economic way. Since the method is converging exponentially with respect to the
tensor rank, the method is also very efficient (but more expensive than the
first two methods for the very first terms in the tensor representation).
However, its implementation requires the realization of inverses of
discretisation matrices in a sparse $\mathcal{H}$-matrix format and a contour
quadrature approximation of the Cauchy-Dunford integral by sinc quadrature by
using a non-trivial parametrisation of the contour.

We expect that these methods can be further developed and an error analysis
which takes into account all error sources (contour quadrature,
discretisation, iteration error, asymptotics with respect to the length of the
cylinder, $\mathcal{H}$-matrix approximation) seems to be feasible. Also the
methods are interesting in the context of a-posteriori error analysis to
estimate the error due to the truncation of the tensor representation at a
cost which is proportional to the solution of problems on the cross sections.
We further expect that more general product domains of the form $%
\BIGOP{\times}%
_{m=1}^{d}\omega_{m}$ for some $\omega_{m}\in\mathbb{R}^{d_{m}}$ with
dimensions $1\leq d_{m}\leq d$ such that $\sum_{m=1}^{d}d_{m}=d$ and domains
with outlets can be handled by our methods since also in this case zero-th
order tensor approximation can be derived by asymptotic analysis (see
\cite{chipot2016_asym}).

\textbf{Acknowledgment}. This work was performed, in part, when the first
author was visiting the USTC in Hefei and during his part time employment at
the S. M. Nikolskii Mathematical Institute of RUDN University, 6
Miklukho-Maklay St, Moscow, 117198. The publication was supported by the
Ministry of Education and Science of the Russian Federation.
\bibliographystyle{abbrv}
\bibliography{mybib_mcwhstasav}

\end{document}